\documentclass[11pt]{article}
\usepackage{epsfig,graphics,graphicx}
\usepackage{amssymb}
\usepackage{amsmath}
\usepackage{extarrows}
\usepackage[utf8]{inputenc}
\usepackage[english]{babel}
\usepackage{geometry}
\usepackage{amsthm}
\usepackage{scalerel}[2014/03/10]
\usepackage[usestackEOL]{stackengine}
\usepackage{esint}
\usepackage{xcolor}

\newtheorem{theo}{Theorem}

\newtheorem{coro}[theo]{Corollary}

\def\R{\mathbb R}

\def\D{\mathbb D}
\def\C{\mathbb C}

\def\K{\mathbb K}
\def\div{\operatorname{div}}

\def\arg{\operatorname{arg}}

\title{\sc{Pointwise rotation \\ for homeomorphisms with integrable distortion \\ and controlled compression}}
\author{Lauri Hitruhin, Banhirup Sengupta}
\date{}

\begin{document}
	
	\maketitle
	
	\abstract{We obtain sharp rotation bounds for homeomorphisms $f:\C\to\C$ whose distortion is in $L^p_{loc}$, $p\geq1$, and whose inverse have controlled modulus of continuity. The motivation to study this class of maps comes from so-called Yudovich solutions to planar Euler equations. Furthermore, we present examples proving sharpness in a strong sense, thereby settling the borderline case $p=1$ in \cite[Theorem 3]{CHS}.} 
	
	\section{Introduction}
	
	An orientation-preserving Sobolev homeomorphism $f\in W^{1,1}_{loc}(\C; \C)$ is called $K$-quasiconformal if the differential $Df(\cdot)$ vanishes in the zero set of the Jacobian determinant $J(\cdot,f)$ and there exists a positive bounded measurable function $\K(\cdot, f):\C\to[1,+\infty]$ such that
	$$|Df(z)|^2\leq K\cdot J(z,f),$$
	at almost every point $z\in\C$. Here,
	$K=\|\K(\cdot, f)\|_\infty$ and $|Df(z)|$ stands for the operator norm of the differential matrix $Df(z)$ of $f$ at the point $z$. If one drops the assumption of $f$ to be a homeomorphism and slightly relax the condition on $\K$ to be only measurable, then $f$ is called a map of finite distortion. Mappings of finite distortion came into existence as a generalization of quasiconformal maps, partially motivated by open questions in nonlinear elasticity. The interested reader is referred to the monograph \cite{AIM} for quasiconformal
	maps in the plane, and to \cite{IM} for a background on mappings of finite distortion.\\
	\\
	Lately there has been a significant upsurge of interest in understanding the pointwise rotational behavior of planar homeomorphisms of finite distortion, along with the spiraling rate of these maps, see \cite{AIPS, B, H1, H2, H3, H4}. To be precise, given such a homeomorphism $f:\C\to \C$ normalized by $f(0)=0$ and $f(1)=1$, one is interested in the maximal growth of $|\arg(f(r))|$ as $r\to 0$. This growth represents the winding number of the image $f([r,1])$ around the origin as $r\to 0^{+}$. It is known that this quantity admits several speeds of growth depending on the class of maps under study. For example, it was shown in \cite{AIPS} that if $f$ is $K$-quasiconformal then
	\begin{equation}\label{qc}
		|\arg(f(r))|\leq \frac12\left(K-\frac1K\right)\,\log\left(\frac1r\right)+c_K,\hspace{1cm}\text{for all }0<r<1.
	\end{equation}
	On the other hand, for homeomorphisms of finite distortion  situation changes drastically and the order of spiraling  depends on integrability of the distortion function. Namely, the first named author discovered in \cite{H2} that if $e^{\K(\cdot, f)}\in L^p_{loc}$ for some $p>0$ then  
	$$|\arg(f(z))|\leq \frac{c}{p}\,\log^2\left(\frac1{|z|}\right),\hspace{1cm}\text{for small enough }|z|,$$
	and moreover this is sharp up to the constant $c>0$. In other words, there is a certain payoff to transit from boundedness to exponential integrability of $\K(\cdot, f)$. Precisely, it is an increment in the exponent of the logarithmic term. Further optimal results were obtained later on in \cite{H3} for homeomorphisms with integrable distortion, that is, when $\K(\cdot, f)\in L^p_{loc}$ for some $p>1$,
	\begin{equation}\label{Lpdist}
		|\arg(f(z))|\leq\frac{c}{|z|^\frac2p},\hspace{1cm}\text{for small enough }|z|
	\end{equation}
	or even when $\K(\cdot, f)\in L^1_{loc}$,
	\begin{equation}\label{L1dist}
		\limsup_{|z|\to0}|z|^{2}|\arg(f(z))|=0.
	\end{equation}
	The above estimates show, as one would expect, that more spiraling can be allowed by relaxing the degree of integrability of $\K(\cdot, f)$.\\
	\\
	One of the most important applications of  mappings of finite distortion is their prominent role in the study of fluid dynamics. 
	Recently it has been proven that there is an important class of homeomorphisms, arising from fluid mechanics, all of which have distortion function in some $L^p_{loc}$ space and are bi-H\"older. This class, first found at \cite{CJ}, contains as an example the flow maps of any Yudovich solution to the planar incompressible Euler system of equations
	\begin{equation}\label{euler}
		\begin{cases}
			\omega_t+ (\mathbf{v} \cdot\nabla )\omega = 0\\
			\div(\mathbf{v})=0\\
			\omega(0,\cdot)=\omega_0.\end{cases}
	\end{equation}
	Here $\omega=\omega(t,z):[0,T]\times\C\to\C$ is the unknown, $\omega_0\in L^\infty(\C;\C)$ is given, and $\mathbf{v}$ is the velocity field. The Biot-Savart law,
	$$\mathbf{v}=\frac{i}{2\pi\overline{z}}\ast\omega$$
	provides the precise relation between $\mathbf{v}$ and $\omega$. The fundamental result by Yudovich \cite{Y} proved existence and uniqueness of a solution $\omega\in L^\infty([0,T]; L^\infty(\C;\C))$ for any given $\omega_0$. In particular, the corresponding velocity field $\mathbf{v}$ belongs to the Zygmund class, and therefore the classical Cauchy-Lipschitz theory guarantees for the ODE
	$$
	\begin{cases}
		\frac{d}{dt}X(t,z)=\mathbf{v}(t, X(t,z))\\
		X(0,z)=z
	\end{cases}
	$$
	both existence and uniqueness of a flow map $X:[0,T]\times\C\to\C$. It was proven in \cite{CJ} that, for small enough $t>0$, each of the flow homeomorphisms $X_t=X(t,\cdot):\C\to\C$ is indeed a mapping of finite distortion. Moreover, for each such small values of $t>0$ there is a number $p(t)>1$ such that the distortion function $\K(\cdot, X_t)$ belongs to $L^p_{loc}$ whenever $p<p(t)$. The H\"older continuity of both $X_t$ and $X_t^{-1}$ was proven in \cite{W}. \\
	\\
	The authors reformulated the bound \eqref{Lpdist} in \cite{CHS} to fit in this situation, getting inspired by this connection between mappings of finite distortion and fluid dynamics . More precisely, it was shown that if $f:\C\to\C$ is a homeomorphism with finite distortion such that $\K(\cdot, f)\in L^p_{loc}$ for some $p>1$ and $f$ has a $\alpha$-H\"older continuous inverse, then
	\begin{equation}\label{chs}
		|\arg(f(z))|\leq C\,\sqrt{\alpha}\, |z|^{-\frac1p}\,\log^\frac12\left(\frac1{|z|}\right),\hspace{1cm}\text{for small enough }|z|,
	\end{equation}
	which drastically improves the bound \eqref{Lpdist} obtained without the H\"older assumption. This is not surprising as the local rotational properties depend not just on the integrability of $\K(\cdot,f)$ but also on the local stretching properties of the mapping, see \cite{AIPS, H2, H3}. In this article, we extend our previous work from \cite{CHS} to a more general class of homeomorphisms, which have $L^p_{loc}$ distortion for $p\geq1$ and the inverse having predetermined modulus of continuity.  \\
	\\
	\begin{theo}\label{main p}
		Let $f:\C\to\C$ be a homeomorphism of finite distortion such that $f(0)=0$, $f(1)=1$, \\and assume that $\K(\cdot, f)\in L^p_{loc};$ $p>1.$ Then
		\begin{equation}\label{p distortion}
			\left|\arg\left(f(z)\right)\right|\leq C\,|z|^{-\frac{1}{p}}\,\log^{\frac{1}{2}}\left(\frac{1}{\underset{\rm |\omega|=|z|}{\rm\min}|f(\omega)|}\right) \quad \text{when $|z|$ is small.}
		\end{equation}
		
		Furthermore, if we assume that $\K(\cdot, f)\in L^{1}_{loc},$ then
		\begin{equation}\label{1 distortion}
			\limsup_{|z|\to0}\frac{|z|}{\sqrt{\log\left(\frac{1}{\underset{\rm |\omega|=|z|}{\rm\min}|f(\omega)|}\right)}}|\arg(f(z))|=0.
		\end{equation}
	\end{theo}
	Towards the optimality of Theorem \ref{main p}, we can show the following.
	\begin{theo}\label{submain p}
		Let $\varphi$ be  a radially increasing homeomorphism with $p$-integrable distortion, $p\geq1$, such that 
		\begin{equation}\label{q condition}
			e^{-g_{\varphi,p}(|z|)|z|^{-\frac{2}{p}}} \leq |\varphi(z)|<|z|^4\quad \text{when $|z|$ is small,} 
		\end{equation}
		where $g_{\varphi,p}:\R \to \R$ is an increasing continuous function with $g(r)\to 0$ when $r \to 0$. Then we can choose an increasing onto homeomorphism $h:[0,+\infty)\to[0,+\infty)$, which can converge to zero as slow as we want, and find a homeomorphism $\overline{f}:\C\to\C$ with the following properties:
		\begin{itemize}
			\item $\overline{f}$ is a homeomorphism of finite distortion, with $\K(\cdot, \overline{f})\in L^p_{loc}$; $\overline{f}(0)=0$, $\overline{f}(1)=1$.
			\item There exists a decreasing sequence $\{r_n\}$,  such that
			\begin{equation}\label{modulus condition}
				|\overline{f}(r_n)|=|\varphi(r_n)|  
			\end{equation}and
			\begin{equation}\label{optimal behavior}
				\left|\arg\left(\overline{f}(r_n)\right)\right|\geq r_{n}^{-\frac{1}{p}}\,\log^{\frac{1}{2}}\left(\frac{1}{|\overline{f}(r_n)|}\right)\,h(r_n).
			\end{equation}
		\end{itemize}
	\end{theo}
	\noindent
	Note that the homeomorphism $\overline{f}$ in Theorem \ref{submain p} is radial and hence $\underset{\rm |\omega|=|z|}{\rm\min}|\overline{f}(\omega)|=|\overline{f}(z)|$. Since $h$ can be chosen to approach zero at any speed, Theorem \ref{submain p} shows that the upper bound provided in Theorem \ref{main p} is essentially sharp when we restrict modulus to satisfy \eqref{q condition}.\\
	\\
	Let us now briefly give some explanation for the bounds \eqref{q condition}. The one on the right specifies that we are studying maps that compress stronger than H\"older maps, which were studied in \cite{CHS}, and thus have faster maximal spiraling rate than given in \eqref{chs}. On the other hand, the bound on the left is always satisfied when $p=1$, see \cite{H3}, and when $p>1$ it is exact up to the gauge function $g_{\varphi, p}$, see \cite{KT}. Studying rotation under extremal compression leads to the extremal pointwise spiraling  as shown in \cite{H3}. Thus Theorem \ref{submain p}, together with examples in \cite{H3} proving optimality of the extremal spiraling rate \eqref{Lpdist}, show that whenever the map $f$ is compressing we have essentially sharp spiraling rates.     \\
	\\
	As a Corollary to Theorem \ref{main p} we can extend the case with Hölder bounds and $p>1$, see  \cite[Theorem 1]{CHS},  to the borderline situation $p=1$.
	
	\begin{coro}\label{main 1}
		Let $f:\C\to\C$ be a homeomorphism of finite distortion such that $f(0)=0$ and $f(1)=1$, and assume that $\K(\cdot, f)\in L^1_{loc}$. Moreover, let us suppose that
		$$|f(x)-f(y)|\geq C\left|x-y\right|^\alpha\hspace{1cm}\text{if }|x-y|\text{ is small,}$$
		for some $\alpha\geq1$. Then
		\begin{equation}\label{Holder 1 distortion}
			\limsup_{|z|\to0}\frac{|z|}{\sqrt{\log\left(\frac{1}{|z|}\right)}}|\arg(f(z))|=0. 
		\end{equation}
	\end{coro}
	Note that in the case $p=1$ we get an improvement in the form of vanishing  limsup compared to the case $p>1$, which is described by the bound \eqref{chs}. This is analogous to the maximal spiraling bounds \eqref{Lpdist} and \eqref{L1dist}, where the exact same improvement happens. \\
	\\
	Finally, we prove the optimality of the above result in a strong sense.
	
	\begin{theo}\label{submain 1}
		Given an increasing, onto homeomorphism $h:[0,+\infty)\to[0,+\infty)$, there exists a homeomorphism $\overline{f}:\C\to\C$ with the following properties:
		\begin{itemize}
			\item[(a)] $\overline{f}$ is a mapping of finite distortion, with $\K(\cdot, \overline{f})\in L^1_{loc}$; $\overline{f}(0)=0$, $\overline{f}(1)=1$.
			\item[(b)] If $\alpha \geq 6,$ then $|\overline{f}(x)-\overline{f}(y)|\geq C|x-y|^\alpha$ whenever $|x-y|<1.$
			\item[(c)] There exists a decreasing sequence $\{r_n\}$, with limit $r_n\to 0^+$ as $n\to \infty$, for which
			$$|\arg(\overline{f}(r_n))|\geq\frac{h(r_n)}{r_n}\, \left(\log\left(\frac1{r_n}\right)\right)^\frac12.$$
		\end{itemize} 
	\end{theo}
	
	\noindent
	Towards the proof of Theorem \ref{submain 1}, we modify the construction from \cite[Theorem 3]{CHS} giving optimality for the bound \eqref{chs}. However, this construction as written in \cite{CHS} does not cover the case $\K\in L^{1}_{loc}$, and thus some changes in the argument are necessary. It turns out that these modifications do not only apply to the $p=1$ setting, and instead work as well when $p>1$. In this case, the H\"older exponent of the inverse map $\overline{f}^{-1}$ from \cite[Theorem 3]{CHS} can be improved from $\frac{p-1}{3p}$ to $\frac{p}{3(p+1)}$. In particular, as $p\searrow 1$ this exponent converges to $\frac{1}{6}$, as one would reasonably expect from Theorem \ref{submain 1}. Also we note that Corollary \ref{main 1} is extremely sharp as homeomorphism $h$ in Theorem \ref{submain 1} can go to zero as slow as we wish. \\
	\\
	\textbf{Acknowledgements}. The authors are partially supported by projects $2017SGR395$ (Govt. of Catalonia), $ PID2020-112881GB-I00$ (Govt. of Spain) and by  the ERC grant 834728 Quamap, the Finnish Academy of Science projects 13316965, 1346562 respectively.

	\section{Preliminaries}
	
	We start by defining pointwise spiraling. Assume $f: \mathbb{C} \to \mathbb{C}$ is a homeomorphism of finite distortion and fix a point $z_0 \in \mathbb{C}$. In order to study the pointwise rotation of $f$ at the point $z_0$ we start by fixing an argument $\theta\in[0,2\pi)$ and study the change in the quantity
	$$\arg (f(z_0+te^{i\theta})-f(z_0))$$ 
	as the parameter $t$ goes from 1 to a small $r$. This can also be understood as the winding of the path $f\left( [z_0+re^{i\theta}, z_0+e^{i\theta}] \right)$ around the point $f(z_0)$. As we are interested in maximal pointwise spiraling, it is natural to find the maximum over all directions $\theta$,
	\begin{equation}\label{spiraling}
		\sup_{\theta \in [0,2\pi)} |\arg (f(z_0+re^{i\theta})-f(z_0)) - \arg (f(z_0+e^{i\theta})-f(z_0))  |.
	\end{equation}
	The maximal pointwise rotation is precisely the behavior of the above quantity \eqref{spiraling} when $r\to 0$. That is, we say that the map $f$ \emph{spirals at the point $z_0$ with a rate $g$}, where $g: [0, \infty) \to [0, \infty)$ is a continuous non-increasing function, if 
	\begin{equation}\label{Spiral rate}
		\limsup_{r\to 0} \frac{\sup_{\theta \in [0,2\pi)} |\arg (f(z_0+re^{i\theta})-f(z_0)) - \arg (f(z_0+e^{i\theta})-f(z_0))  |}{g(r)} = C
	\end{equation}
	for some constant $0<C<\infty$. Finding maximal pointwise rotation for a given class of maps is equivalent to find maximal growth of function $g$ when $r \to 0$  for this class. \\
	\\
	The proof of Theorem \ref{main p} relies heavily on the modulus of path families. We provide the main definitions, see, for example, \cite{V} for a closer look at the topic. An image of a line segment $I$ under a continuous map is called a \emph{path}, and we denote by $\Gamma$  a family of paths. Given a path family $\Gamma$, we say that a Borel measurable function $\rho$ is \emph{admissible for $\Gamma$} if any rectifiable $\gamma \in \Gamma$ satisfies
	\begin{equation*}
		\int_{\gamma} \rho(z)dz \geq 1.
	\end{equation*}
	The \emph{modulus of the path family $\Gamma$} is then defined by 
	\begin{equation*}
		M(\Gamma)= \inf_{\rho \text{ admissible}} \int_{\mathbb{C}} \rho^{2}(z)\,dA(z).
	\end{equation*}
	We will also need a weighted version of the modulus. Any measurable, locally integrable function $\omega: \mathbb{C} \to [0, \infty)$ will be called \emph{weight function}. In our case the weight function $\omega$ will always be the distortion function $\K(\cdot, f)$. Then, we define the weighted modulus $M_{\omega}(\Gamma)$ by 
	\begin{equation*}
		M_{\omega}(\Gamma)= \inf_{\rho \text{ admissible}} \int_{\mathbb{C}} \rho^{2}(z)\, \omega(z)\,dA(z).
	\end{equation*}
	Finally, we need the modulus inequality
	\begin{equation}\label{Modeq}
		M(f(\Gamma)) \leq M_{\K(\cdot, f)}(\Gamma)
	\end{equation}
	which holds for any homeomorphism $f$ of finite distortion for which the distortion $\K(\cdot, f)$ is locally integrable, proven by the first named author in \cite{H3}.
	
	\section{Spiraling bounds}
	In this section, we prove Theorem \ref{main p} following closely proof for extremal spiraling rate \eqref{chs} in \cite{CHS}.
	\noindent
	\begin{proof}
		Let $f$ satisfy the hypothesis of Theorem \ref{main p}, and  let $z\in\C\setminus\{0\}$ be such that $|z|<1.$ Our goal is to estimate the \emph{winding number} $n(z)$ of the image set $f\left(\left[z,\frac{z}{|z|}\right]\right)$ around the origin (recall that $f(0)=0$). We will bound $n(z)$ using the modulus inequality \eqref{Modeq}. More precisely, we will prove that  
		\begin{equation}
			n(z)\leq C\,|z|^{-\frac{1}{p}}\,\log^{\frac{1}{2}}\left(\frac{1}{\underset{\rm |z_0|=|z|}{\rm\min}|f(z_0)|}\right),
			\text{ when $p>1$ and } \limsup_{|z|\to0}\frac{|z|}{\sqrt{\log\left(\frac{1}{\underset{\rm |z_0|=|z|}{\rm\min}|f(z_0)|}\right)}}n(z)=0 \text{ when $p=1$. } 
		\end{equation}
		Let us first prove $p>1$ case. To this end, choose an arbitrary point $z_0\in\C\setminus\{0\}$ such that $|z_0|<1$. Without loss of generality we can assume that $z_0$ is real and positive. Next, fix line segments $E=[z_0,1]$ and $F=(-\infty,0]$, and let $\Gamma$ be the family of paths connecting them. Then, following the proof of Theorem 1 in \cite{CHS} (page numbers 5 and 6), we can similarly estimate the modulus term $M_{\K(\cdot, f)}(\Gamma)$ from above as, 
		\begin{equation}\label{modulus with distortion p}
			M_{\K(\cdot, f)}(\Gamma)\leq c_{f,p}z_{0}^{-\frac{2}{p}}.
		\end{equation}
		Next, we would like to estimate the modulus term $M\left(f(\Gamma)\right)$ from below for $p\geq1$. Let us recall that $f(0)=0$ and write $M\left(f(\Gamma)\right)$ in polar coordinates as follows:
		$$\aligned
		M\left(f(\Gamma)\right)&=\inf_{\rho\hspace{0.1cm}admissible}\int_{\C}\rho^{2}(z)\,dA(z)
		=\inf_{\rho\hspace{0.1cm}admissible}\int_{0}^{2\pi}\int_{0}^{\infty}\rho^{2}(r,\theta)r\hspace{0.1cm}drd\theta.\\
		\endaligned$$ Let us provide a lower bound for 
		$$\aligned
		\int_{0}^{\infty}\rho^{2}(r,\theta)rdr
		\endaligned$$ for an arbitrary direction $\theta\in[0,2\pi)$ and an arbitrary admissible $\rho$. It is known from the proof of Theorem 1 in \cite{CHS} (page numbers 6 and 7), that
		\begin{equation}\label{Lower bound}
			\int_{0}^{\infty}\rho^{2}(r,\theta)rdr\geq\frac{n^{2}(z_0)}{\log\left(\frac{c_f}{r_f}\right)},
		\end{equation}
		where $c_f=\sup_{z\in E}|f(z)|$ and $r_f=\inf_{z\in E}|f(z)|$.
		The constant $c_f$ is finite and does not depend on either $\theta$ or $z_0$, at least for small $z_0$. So, it is irrelevant at the limit $z_0\to0$. Hence the estimate \eqref{Lower bound} implies that
		\begin{equation}\label{Lower bound general modulus}
			M\left(f(\Gamma)\right)\geq\frac{c n^{2}(z_0)}{\log\left(\frac{1}{\underset{\rm |z|=|z_0|}{\rm\min}|f(z)|}\right)}.
		\end{equation}
		Next, use the modulus inequality \eqref{Modeq} to get
		$$\aligned
		\frac{n^{2}(z_0)}{\log\left(\frac{1}{\underset{\rm |z|=|z_0|}{\rm\min}|f(z)|}\right)}\leq c_{f,p}z_{0}^{-\frac{2}{p}},
		\endaligned$$which implies the desired estimate \eqref{p distortion}.\\
		\noindent
		\textbf{To prove $p=1$ case}, we will again use the modulus inequality \eqref{Modeq}. Note that we have already lower bound for $M\left(f(\Gamma)\right)$ from \eqref{Lower bound general modulus} for any $p\geq1$. Therefore, we just need to estimate modulus term $M_{\K(\cdot, f)}(\Gamma)$ from above. To this end, let us define the function
		$$\aligned
		\rho_{0}(z)=
		\begin{cases}
			\frac{1}{z_0}\hspace{0.4cm}\mbox{if}\hspace{0.2cm}dist(z,E)<z_0\\
			0\hspace{0.6cm}otherwise\\
		\end{cases}
		\endaligned$$
		Note that $\rho_0$ is admissible with respect to the path family $\Gamma$. Therefore,
		$$\aligned
		M_{\K(.,f)}(\Gamma)&\leq\int_{\C}\K(.,f)\rho_{0}^{2}(z)dA(z)
		=\frac{1}{z_{0}^{2}}\int_{\{z:dist(z,E)<z_0\}}\K(.,f)(z)dA(z).\\
		\endaligned$$
		Denote $$\int_{\{z:dist(z,E)<z_0\}}\K(.,f)(z)dA(z)= C_f(z_0)$$ and note that since $\K(.,f)(z)\in L^{1}_{loc}(\C)$ and $\left|\{z:dist(z,E)<z_0\}\right|\to0$ it follows that $C_{f}(z_0)\to0$ as $z_0\to0$, and thus
		\begin{equation}\label{modulus with distortion 1}
			M_{\K(.,f)}(\Gamma)\leq\frac{C_{f}(z_0)}{z_{0}^{2}}.
		\end{equation}
		Next, we use the modulus inequality \eqref{Modeq}, bounds \eqref{Lower bound general modulus} and \eqref{modulus with distortion 1} to get
		$$\aligned
		\frac{n^{2}(z_0)}{\log\left(\frac{1}{\underset{\rm |z|=|z_0|}{\rm\min}|f(z)|}\right)}\leq\frac{C_{f}(z_0)}{z_{0}^{2}}
		\endaligned$$which implies the desired estimate \eqref{1 distortion}. Hence, Theorem \ref{main p} is proved.
	\end{proof}
	\noindent
	Corollary \ref{main 1} follows directly from Theorem \ref{main p} by estimating $\underset{|\omega|=|z|}{\min}|f(\omega)|$ from below by H\"older continuity.
	
	\section{Optimal Results}
	We are going to prove Theorems \ref{submain p} and \ref{submain 1} in this section.
	\noindent
	\begin{proof}[Proof of Theorem \ref{submain p}]
		
		As a first step we construct a map which \emph{only rotates}. This map will have the correct spiraling rate but the distortion of the map will not belong to the desired space. To overcome this barrier we compose it with a radial stretching map, which gives us better control over distortion.\\
		\\
		Given an arbitrary annulus $A=B(0,R)\setminus B(0,r)$ let us define the corresponding rotation map as 
		\begin{equation}\label{rotation building block}
			\phi_{A}(z)=
			\begin{cases}
				z&|z|>R\\
				z \,e^{ i\alpha\log\left|\frac{z}{R}\right| }&r\leq |z|\leq R\\
				z\,e^{i\alpha \,\log\frac{r}{R} }&|z|<r
			\end{cases}
		\end{equation}
		Here $0<r<R$, and $\alpha \in\R$. It is clear that $\phi_{A}:\C\to\C$ is bilipschitz, hence quasiconformal (its quasiconformality constant depends only on $\alpha$), and moreover it is conformal outside the annulus $A$. Moreover, $|\phi_A(te^{i\theta})|=t$ for each $t>0$ and $\theta\in\R$. This means that $\phi_A$ leaves fixed all circles centered at $0$. It is easy to check that the jacobian determinant $J(z, \phi_A)=1$ for each $z$. \\
		\\
		Next, let us consider sequence $\{r_n\}$ such that $0<r_{n+1}<\frac{r_n}{2e}$ and $r_1<\frac{1}{e}$. Also, let $R_n=e r_n$, which ensures that $2r_{n+1}<R_{n+1}<\frac{r_n}2$. Let us now construct disjoint annuli $A_n=B(0,R_n)\setminus B(0,r_n)$, and set $\{f_n\}_n$ to be a sequence of maps, constructed in an iterative way as follows. For $n=1$, we set
		\begin{equation}\label{initial rotation}
			f_{1}(z)=\phi_{A_1}(z)=
			\begin{cases}
				z&|z|>R_1\\
				z \,e^{ i\alpha_1\,\log \frac{|z| }{R_1}}&r_1\leq|z|\leq R_1\\
				z\,e^{-i\alpha_1 }&|z|<r_1\\
			\end{cases}
		\end{equation}
		where $\alpha_1\in\R$, $\alpha_1\geq 1$, is to be determined later. We then define $f_n$ for $n\geq2$ as 
		$$\aligned
		f_{n}(z)=\phi_{f_{n-1}(A_n)}\circ f_{n-1}(z)
		\endaligned$$
		again for some values $\alpha_n\in\R$, $\alpha_n\geq 1$, to be determined later. Clearly, each $f_n:\C\to\C$ is quasiconformal, and conformal outside the annuli $A_i$, $i=1,\cdots, n$. It is also clear that $f_n(z)=f_{n-1}(z)$ on the unbounded component of $\C\setminus f_{n-1}(A_n)$ (i.e. outside of $B(0, R_n)$). This proves that the sequence $f_n$ is uniformly Cauchy and hence it converges to a map $f$, which is again a homeomorphism by construction. Now, since $f_n$ is quasiconformal for every $n$ and $f_n(z)=f_{n-1}(z)$ everywhere except inside the ball $B(0,R_n)$, where $R_n\to0$ as $n\to\infty$, the limit map $f$ is absolutely continuous on almost every line parallel to the coordinate axes and differentiable almost everywhere. Direct calculation shows that
		$$
		|D\phi_{A_n}(z)|=|\partial\phi_{A_n}(z)|+|\overline\partial\phi_{A_n}(z)|=
		\begin{cases}
			1&|z|>R_n\\
			\frac{|2+i\alpha_n|+|\alpha_n|}{2}&r_n\leq |z|\leq R_n\\
			1&|z|< r_n\end{cases}
		$$
		which allows us to estimate that
		$$\aligned
		|\partial f(z)|+|\overline\partial f(z)|\leq 2\alpha_n\hspace{1cm}\text{whenever }z\in A_n,
		\endaligned$$
		and $|Df(z)|\leq1$ otherwise. Therefore, in order to have $Df(z)\in L^{1}_{loc}(\C)$ it suffices that
		\begin{equation}\label{dfinL1}
			\sum_n \alpha_n\,r_n^2<+\infty.
		\end{equation}
		This, together with the absolute continuity, guarantees $f\in W^{1,1}_{loc}(\C)$. Also, since $f$ is a homeomorphism, we have that $J(\cdot,f)\in L^{1}_{loc}(\C)$, and in fact $J(z, f)=1$ at almost every $z\in\C$. Therefore, $f$ is a homeomorphism of finite distortion, with distortion function
		$$\aligned
		\K(z,f)=\frac{|Df(z)|^2}{J(z, f)}\leq
		\begin{cases}
			4\alpha_n^2&z\in A_n,\\
			1&\text{otherwise.}
		\end{cases}
		\endaligned$$
		Especially, in order to have $\K(\cdot, f)\in L^p_{loc}$, it suffices to ensure the convergence of the series
		\begin{equation}\label{kinLp}
			\sum_{n=1}^{\infty}|A_n|(4\alpha_n^2)^{p}\simeq \sum_{n=1}^\infty \alpha_n^{2p}\,r_n^2
		\end{equation}
		which can be done by choosing $\alpha_n$ properly. Note that if \eqref{kinLp} holds, then also \eqref{dfinL1} holds, because our choice of $\alpha_n$ will guarantee $\alpha_n\geq 1$. The last restriction to choose our $\alpha_n$ comes from rotational behavior of $f$. It is clear from the above construction that $f(0)=0$, $f(1)=1$ and
		$$
		\left|\arg\left(f(r_n)\right)\right|\geq\left|\arg\left(\left(\frac{1}{e}\right)^{1+i\alpha_n}\right)\right|=\alpha_n
		$$
		for every $r_n$. Let us choose $\alpha_n= r_n^{-1/p}\,\log^{1/2}(1/\varphi(r_n))$. This implies that 
		$$\aligned
		\left|\arg\left(f(r_n)\right)\right|\geq r_n^{-1/p}\,\log^{1/2}(1/\varphi(r_n))
		\endaligned,$$ which shows that this map would be optimal for Theorem \ref{main p}. However, with this particular choice of $\alpha_n$, 
		$$\aligned
		\sum_{n=1}^\infty\, \alpha_n^{2p}\,r_n^2=\sum_{n=1}^\infty\, \log^{p}(1/\varphi(r_n))
		\endaligned$$ which is certainly not finite. This means that $\K(\cdot, f)\notin L^p_{loc}$. \\
		\\
		Hence we need to modify the construction by adding a stretching factor to our building blocks, which lets us reduce the distortion while preserving spiraling rate. This is precisely done by substituting the logarithmic spiral map $z|z|^{i\alpha}=ze^{i\alpha\log|z|}$ by a complex power $z|z|^{q+i\alpha}=z|z|^q\,e^{i\alpha\log|z|}$ at each iterate. Similarly as in the previous construction, we consider a rapidly decreasing sequence $\{r_n\}$ such that $r_{n+1}<\frac{r_n}{2e}$, $r_1<\frac{1}{e}$ and set $R_n=e r_n$. Given an arbitrary annulus $A=B(0,R)\setminus B(0,r)$ we define the corresponding composition map as follows:
		\begin{equation}\label{buildingblock}
			\aligned
			\phi_{A}(z)=
			\begin{cases}
				z&|z|>R\\
				z\,\left|\frac{z}{R}\right|^{q-1}\,e^{i\alpha\log\frac{|z|}{R}}&r\leq|z|\leq R\\
				z\left(\frac{r}{R}\right)^{q-1}\,e^{i\alpha\log\frac{r}{R}}&|z|<r\
			\end{cases}
			\endaligned
		\end{equation} 
		Note that we will always choose $q\geq 1$. Direct calculation shows that
		\begin{equation}\label{differential}
			|\partial\phi_A(z)|+|\overline\partial\phi_A(z)|
			=\begin{cases}
				1&|z|>R\\
				R^{1-q}|z|^{q-1}\frac{|q+1+i\alpha|+|q-1+i\alpha|}{2}&r \leq |z|\leq R \\
				R^{1-q}r^{q-1}&|z|<r
			\end{cases} 
		\end{equation}
		\noindent
		and 
		\begin{equation}\label{distortion}
			\K(z,\phi_A)
			=\begin{cases}
				1&|z|>R\\
				\frac{(|q+1+i\alpha|+|q-1+i\alpha|)^2}{4q}&r \leq |z|\leq R \\
				1&|z|<r
			\end{cases}    
		\end{equation}
		In particular, if $2\leq q+1 \leq \alpha$, which will be satisfied for our choices of $\alpha$ and $q$, then one may estimate $\|\K(\cdot, \phi_A)\|_\infty\leq \frac{4\alpha^2}{q}$. Next, let us construct the sequence of maps $f_n$ in an iterative way. For $n=1$, we set
		\begin{equation}\label{initial rotstretch}
			f_{1}(z)=\phi_{A_1}(z)=
			\begin{cases}
				z&|z|<R_1\\
				z\,\left|\frac{z}{R_1}\right|^{q_1-1}\,e^{i\alpha_1\log\frac{|z|}{R_1}}&r_1\leq |z|\leq R_1\\
				z\left(\frac{1}{e}\right)^{q_1-1}\,e^{-i\alpha_1}&|z|<r_1\\
			\end{cases}    
		\end{equation}
		where $q_1$ and $\alpha_1$ are to be determined later. Assuming we have $f_1,\cdots, f_{n-1}$, we define $f_n$($n\geq2$) as:
		$$\aligned
		f_{n}(z)=\phi_{f_{n-1}(A_n)}\circ f_{n-1}(z).
		\endaligned$$
		Similarly as in the previous construction each $f_n:\C\to\C$ is quasiconformal, and conformal outside the annuli $A_i$, $i\in\{1,...,n\}$. Moreover, one can easily show that 
		\begin{equation}\label{distcomp}
			\K(\cdot , f_n) = \prod_{j=1}^{n}\K(\cdot,f_{n-j} \circ\phi_{f_{n-j}(A_{n-j+1})})\\
			=\prod_{j=1}^{n}\K( \cdot,\phi_{ A_{n-j+1} })    
		\end{equation}
		so that $\K(z, f_n)\leq C\frac{\alpha_j^2}{q_j}$ whenever $z\in A_j$, $j=1,..., n$ while $\K(\cdot,f_n)=1$ otherwise. In a similar way, we can use that $|D\phi_A(z)|\leq C\alpha $ when $z\in A$ (and $|D\phi_A(z)|\leq 1$ at all other points) to obtain that $|Df_n|\leq C \alpha_j$ on $A_j$, $j=1,..., n$, and $|Df_n|\leq 1$ otherwise. By construction, we have $f_n(z)=f_{n-1}(z)$ whenever $z\notin B(0, R_n)$. Thus $\{f_n \}_n$ converges uniformly to a map $\overline{f}$, that is,
		$$\aligned
		\overline{f}=\lim_{n\to\infty}f_n
		\endaligned$$
		which is again a homeomorphism by construction. A similar argument to the one before shows that $\overline{f}$ is absolutely continuous on almost every line parallel to the coordinate axis. For almost every fixed $z_0$ there is a neighborhood of $z_0$ such that the sequence $\{f_n(z)\}_n$ remains constant for $n$ very large and $z$ in that neighborhood. Therefore the same happens to the sequences $Df_n(z)$, $J(z, f_n)$ and $\K(z, f_n)$, and so their limits are precisely $D\overline{f}(z)$, $J(z, \overline{f})$ and $\K(z, \overline{f})$. Especially, in order to have $D\overline{f}\in L^1_{loc}$ it suffices that
		\begin{equation}\label{dfinL1general}
			\sum_{n=1}^\infty |A_n|\,\alpha_n<+\infty.
		\end{equation}
		In case this holds true, then $\overline{f}$ is a homeomorphism in $W^{1,1}_{loc}$, and as a consequence its jacobian determinant $J(\cdot, \overline{f})\in L^1_{loc}$. Moreover, in order to have $\K(\cdot, \overline{f})\in L^p_{loc}$; ($p\geq1$) one needs to require that
		\begin{equation}\label{kfinLpgeneral}
			\sum_{n=1}^\infty |A_n|\,\frac{\alpha_n^{2p}}{q_n^p}<+\infty.
		\end{equation}
		Again, as it was the case for the pure rotation example, when $p>1$ condition \eqref{kfinLpgeneral} implies \eqref{dfinL1general} if $q_n^\frac{p}{2p-1}\leq \alpha_n$, and for $p=1$ case we must verify $q_n\leq\alpha_n$. So, our parameters $\alpha_n$ and $q_n$ need to be chosen according to these constrains as well as the purpose of $\overline{f}$ to be optimal for Theorem \ref{main p}. To this end, note that $\overline{f}(0)=0$, $\overline{f}(1)=1$ and
		\begin{equation}\label{optimality}
			\left|\arg\left(\overline{f}(r_n)\right)\right|\geq\,\left|\arg\left(\left(\frac{1}{e}\right)^{q_n+i\alpha_n}\right)\right|=|\alpha_n|,
		\end{equation}
		which motivates us to choose
		\begin{equation}\label{Precise q_n}
			q_n=
			\begin{cases}
				\log\frac{e\,r_1}{|\varphi(r_1)|}& n=1\\
				\log\left(\frac{e\,\cdot\,r_{n}\cdot\,\left(\frac{1}{e}\right)^{q_{n-1}+q_{n-2}+...+q_1-(n-1)}}{|\varphi(r_n)|}\right)& n\geq2
			\end{cases}
		\end{equation}and
		\begin{equation}\label{Precise a_n}
			\alpha_n= h(r_n)\,\left(\log\frac{1}{|\varphi(r_n)|}\right)^{1/2}\,r_n^{-\frac{1}{p}}   
		\end{equation}
		where $h$ is a monotone non-increasing gauge function such that $h(r)\to0$ as $r\to0$ which we specify later. Next, we show that $q_n\leq\alpha_n$ for all $p \geq 1$, from which $q_n^{\frac{p}{2p-1}}\leq\alpha_n$ also follows. At this point, we impose an ansatz on $r_n$:
		\begin{equation}\label{Additional assumption Lauri}
			r_n<\left(\frac{1}{e}\right)^{q_{n-1}+q_{n-2}+...+q_{1}-(n-1)},    
		\end{equation}
		which is feasible as the radii $r_n$ can be assumed to decrease as fast as we want. Let us recall our assumption on $\varphi$ to satisfy compression bound:
		\begin{equation}\label{CompressionaddedbyLauri}
			|\varphi(z)|\geq e^{-g_{\varphi,p}(|z|)|z|^{-\frac{2}{p}}}, 
		\end{equation}
		where $g_{\varphi,p}:\R\to\R$ is some increasing gauge function such that $|g_{\varphi,p}|\to0$ as $|z|\to0$.
		We calculate 
		$$\aligned
		q_n &=\log\,\left(\frac{e\,\cdot r_n\,\cdot\left(\frac{1}{e}\right)^{q_{n-1}+....+q_1-(n-1)}}{|\varphi(r_n)|}\right)\\
		&\leq\log\,\frac{1}{|\varphi(r_n)|}\\
		&=\log^{\frac{1}{2}}\frac{1}{|\varphi(r_n)|}\cdot\,\log^{\frac{1}{2}}\frac{1}{|\varphi(r_n)|}\\
		&\leq\log^{\frac{1}{2}}\frac{1}{|\varphi(r_n)|}\cdot\, \sqrt{g_{\varphi,p}(r_n)}\frac{1}{r_n^{\frac{1}{p}}}\\
		&\leq\log^{\frac{1}{2}}\frac{1}{|\varphi(r_n)|}\cdot\,\frac{h(r_n)}{r_n^{\frac{1}{p}}}= \alpha_n
		\endaligned$$
		where the last inequality holds for $h$ converging to zero slowly enough. Note that, from \cite[Theorem 1.6]{H3} we see that if $p=1$ then the compression bound \eqref{CompressionaddedbyLauri} is always satisfied with some $g_{\varphi}$. Thus our choices for $q_n$ and $\alpha_n$ satisfy technical constrains. Next, we show that estimate \eqref{kfinLpgeneral} governing integrability of the distortion holds true for $p\geq1$. We start by estimating
		$$\aligned
		|A_n|\frac{\alpha_n^{2p}}{q_n^p}&=C\,r_n^2\,\frac{h^{2p}\,r_n^{-2}\,\log^{p}\left(\frac{1}{|\varphi(r_n)|}\right)}{\log^p\,\left(\frac{e\,\cdot r_n\,\cdot\left(\frac{1}{e}\right)^{q_{n-1}+....+q_1-(n-1)}}{|\varphi(r_n)|}\right)}
		\leq C\,h^{2p}\,\frac{\log^{p}\left(\frac{1}{|\varphi(r_n)|}\right)}{\log^p\,\left(\frac{r_n^2}{|\varphi(r_n)|}\right)}.
		\endaligned$$
		Using the condition \eqref{q condition} we see, $1<\frac{\log^{p}\left(\frac{1}{|\varphi(r_n)|}\right)}{\log^p\,\left(\frac{r_n^2}{|\varphi(r_n)|}\right)}\leq 2^p$, and therefore condition \eqref{kfinLpgeneral} is equivalent to 
		$$\sum_nh(r_n)^{2p}<+\infty$$
		which we can always satisfy by choosing $r_n$ small enough. Having \eqref{kfinLpgeneral} fulfilled, our map $\overline{f}$ is a mapping of finite distortion with $\K(\cdot, \overline{f})\in L^p_{loc}$. Next we must show that our mapping $f$ has right compression and spiraling behavior. Let us start with modulus and show that $|\overline{f}(r_n)|=|\varphi(r_n)|$ by calculating
		$$\aligned
		|\overline{f}(r_n)|&=\left(\frac{1}{e}\right)^{q_{n-1}+...+q_1-(n-1)}\cdot r_n\,\left(\frac{r_n}{R_n}\right)^{q_n-1}
		&=\left(\frac{1}{e}\right)^{q_{n-1}+...+q_1-(n-1)}\cdot r_n\,\left(\frac{1}{e}\right)^{q_n-1}
		=|\varphi(r_n)|,
		\endaligned$$ 
		where the last line follows from the penultimate due to our choice of $q_n$.\\
		\\
		To finish the proof we must show that the rotation bound \eqref{optimal behavior} holds true. But this follows directly from \eqref{optimality}, \eqref{Precise a_n} and from the above modulus equation. Finally we note that in the proof we only need to assure that the sequence $r_n$ decreases fast enough, and thus we can choose it to be a subsequence for an arbitrary predefined sequence $\{\lambda_n\}$. This concludes the proof of Theorem \ref{submain p}. 
	\end{proof}
	\begin{proof}[Proof of Theorem \ref{submain 1}]
		We  prove Theorem \ref{submain 1} in two steps similarly to Theorem \ref{submain p}. In the first step we construct a map which \emph{only rotates}. This map already provides the optimal result in the exponent scale. Then, as in previous construction, we compose this map with radial stretching map and finish the proof.\\
		\\
		Given an arbitrary annulus $A=B(0,R)\setminus B(0,r)$ we define corresponding rotation map $\phi_{A}$ as in \eqref{rotation building block}. It is clear that $\phi_{A}:\C\to\C$ is quasiconformal, and moreover  it is conformal outside the annulus $A$. Furthermore, $\phi_A$ leaves fixed all circles centered at $0$, and Jacobian determinant $J(z, \phi_A)=1,$ for each $z$. Next, consider again sequence $\{r_n\}$ such that $0<r_{n+1}<\frac{r_n}{2e}$, $r_1<\frac{1}{e}$, and fix $R_n=e r_n$. We then construct disjoint annuli $A_n=B(0,R_n)\setminus B(0,r_n)$ and a sequence of maps $\{f_n\}_n$ iteratively as before. That is, set $f_1$ as in \eqref{initial rotation} and 
		define $f_n$ for $n\geq2$ as 
		$$\aligned
		f_{n}(z)=\phi_{f_{n-1}(A_n)}\circ f_{n-1}(z)
		\endaligned$$
		for some $\alpha_n\geq 1$, to be determined later. One can use exactly same arguments as in the proof of Theorem \ref{submain p} to deduce that the limit
		$
		f=\lim_{n\to\infty}f_n
		$
		is a homeomorphism with integrable distortion if 
		\begin{equation}\label{dfinl1}
			\sum_n \alpha_n\,r_n^2<+\infty
		\end{equation}
		and
		\begin{equation}\label{kinl1}
			\sum_{n=1}^{\infty}|A_n|4\alpha_n^2\simeq \sum_{n=1}^\infty \alpha_n^{2}\,r_n^2<+\infty.
		\end{equation}
		In fact, \eqref{kinl1} implies \eqref{dfinl1} as we will choose $\alpha_n\geq1.$ Hence it is sufficient to choose $\alpha_n$ so that \eqref{kinl1} is satisfied. 
		Furthermore, it is clear from above construction that $f(0)=0$, $f(1)=1$ and
		\begin{equation}\label{LauriADD}
			\left|\arg\left(f(r_n)\right)\right|\geq\left|\arg\left(\left(\frac{1}{e}\right)^{1+i\alpha_n}\right)\right|=\alpha_n    
		\end{equation}
		for every $r_n$. Since we want our map to be optimal for Corollary \ref{main 1}, we may be tempted to choose $\alpha_n= \frac{\log^{1/2}(1/r_n)}{r_n}$. Unfortunately such a choice does not meet the requirement \eqref{kinl1} and instead one is forced to choose
		$\alpha_n=\frac{h(r_n)}{r_n},$
		where $h:[0,\infty)\to[0,\infty)$ is a monotonically decreasing gauge function such that $\lim_{r\to 0^+}h(r)=0$.  With this choice, \eqref{kinl1} is fulfilled if
		$$\sum_{n=1}^\infty \left(h(r_n)\right)^{2}<+\infty,$$
		which we can ensure by choosing small enough $r_n$. Note that this does not provide optimality for Corollary \ref{main 1} in full generality, but it already gives the right order in the exponent scale.\\
		\\
		Finally, we show that $f^{-1}$ is H\"older continuous with exponent $\frac{1}{2}$. To this end, let us recall that our map $f$ is actually a limit of iterates of logarithmic spiral maps inside the annuli $A_n=B(0,R_n)\setminus B(0,r_n)$. In particular, as shown in \cite{AIPS}, if $\gamma\in\R$ then the basic logarithmic spiral map $g(z)=ze^{i\gamma\log|z|}$ is $L$-bilipschitz for a constant $L$ such that $|\gamma|=L-\frac{1}{L}$. And thus for large $|\gamma|$ one roughly has $|\gamma|\simeq L$. Since our $f_n$ behaves in the annulus $A_n$ as a spiral map with $|\gamma|=\alpha_n$, we deduce that the bilipschitz constant of $f_n$ on $A_n$ is 
		$L \simeq|\gamma| =\alpha_n= \frac{h(r_n)}{r_n}.$\\ 
		\\
		Let us now start the proof. We first consider the case where $x, y\in A_n$, and hence $f(x)=f_n(x)$, $f(y)=f_n(y)$. Since $r_n>C|x-y|$, we have
		$$\aligned
		|f(x)-f(y)|=|f_n(x)-f_n(y)| &\gtrsim \frac{r_n}{h(r_n)}|x-y|
		\geq \frac{C}{h(r_n)}|x-y|^2
		\geq C|x-y|^2
		\endaligned$$ 
		where we have used the bilipschitz nature of $f_n$ on $A_n$. The fact that $f$ is H\"older from below inside the annuli $A_n$ with exponent $2$ implies that in these sets $f^{-1}$  is H\"older continuous with exponent $\frac{1}{2}.$ Since $f$ and $f^{-1}$ are essentially the same mapping modulo the direction of rotation, $f$ is also $\frac12$ H\"older continuous inside $A_n.$ Then we assume that $x, y\in D_n=B(0,r_n)\setminus B(0,R_{n+1})$. In this case $f$ is of the form $ze^{i\beta}$, where $\beta\in\R\setminus\{0\}$, which is clearly an isometry and hence H\"older estimate inside $D_n$ is trivial.  \\
		\\
		Next, we take $x\in A_n$ and $y\in D_n$. In particular, $|x|\geq|y|$. Then let $w$ be the point on the outer boundary of $D_n$ joining $x$ and $y$. We have
		$$\aligned
		|f(x)-f(y)|\leq|f(x)-f(w)|+|f(w)-f(y)|
		\leq C|x-w|^{\frac{1}{2}}+|w-y|
		\leq2C|x-y|^\frac{1}{2}.
		\endaligned$$
		The same happens if $x\in D_{n-1}$ and $y\in A_n$. So it just remains to see what happens when points are further apart from each other. Let us first cover the case $x\in A_n=B(0,R_n)\setminus B(0,r_n)$ and $y\in B(0,R_{n+1})$. Let $L$ be the line joining $x$ and $y$. We divide it into three parts, viz., $L_1$, $L_2$ and $L_3$. Fix $L_1$ so that it connects $x$ to a point $a$ on the inner boundary of $A_n$, giving estimate
		$$|f(x)-f(a)|=|f_n(x)-f_n(a)|\leq C|x-a|^{\frac12}.$$
		Next, $L_2$ connects $a$ to the crossing point of the line $L$ and the inner boundary of $D_n$, which we denote by $b$. And since $f$ is an isometry in $D_n$ an estimate for line segment $L_2$ is trivial. For $L_3$ part we note that from  $2R_{n+1}<r_n<\frac{R_n}2$ we get that $|f(a)|>2|f(b)|$ and hence
		$$|f(b)-f(y)|\leq2|f(b)|\leq2|f(b)-f(a)|=2|b-a|.$$ 
		Combining these estimates we get
		$$\aligned
		|f(x)-f(y)|&\leq|f(x)-f(a)|+|f(a)-f(b)|+|f(b)-f(y)|
		\leq C|x-a|^{\frac{1}{2}}+|a-b|+2|a-b|
		\leq C|x-y|^{\frac{1}{2}}.
		\endaligned$$
		The case $x\in D_n$ and $y\in B(0,r_{n+1})$ can be proved similarly. Thus $f$ is H\"older continuous with exponent $\frac12$. Here we again note that $f$ and $f^{-1}$ are essentially the same mapping modulo the direction of rotation, and hence $f^{-1} $ is also H\"older continuous with exponent $\frac{1}{2}$. Thus $f$ is H\"older from below  with exponent $2.$ \\
		\\
		As we discussed before, the above example approaches the borderline stated in Corollary \ref{main 1}, but it does not attain full optimality yet. To this end, we need to modify it by adding a stretching factor to our building blocks, which lets us increase rotation without increasing the distortion. This is done by replacing, at each iterate, the logarithmic spiral map $z|z|^{i\alpha}=ze^{i\alpha\log|z|}$ by a complex power $z|z|^{q+i\alpha}=z|z|^q\,e^{i\alpha\log|z|}$. So, similarly as in the previous construction, we consider a rapidly decreasing sequence $\{r_n\}$ such that $r_{n+1}<\frac{r_n}{2e}$, $r_1<\frac{1}{e}$ and fix $R_n=e r_n$. Given an arbitrary annulus $A=B(0,R)\setminus B(0,r)$ we define the corresponding radial stretching combined with rotation map as in \eqref{buildingblock}. As before we will choose $q\geq 1$. \\
		\\
		The values of the differential and distortion of $\phi_{A}$ are already known from \eqref{differential} and \eqref{distortion}.
		In particular, if $2\leq q+1<\alpha$ then one may estimate $\|\K(\cdot, \phi_A)\|_\infty\leq \frac{4\alpha^2}{q}$. Next, we construct the sequence of maps $f_n$ in an iterative way as before. Let us set $f_1$ as in \eqref{initial rotstretch} and $f_n$ for $n\geq2$ as:
		$$\aligned
		f_{n}(z)=\phi_{f_{n-1}(A_n)}\circ f_{n-1}(z).
		\endaligned$$
		Each $f_n:\C\to\C$ is quasiconformal, and conformal outside the annuli $A_i$, $i\in\{1,...,n\}$. Moreover, we still calculate distortion as in \eqref{distcomp}, so that $\K(z, f_n)\leq C\frac{\alpha_j^2}{q_j}$ whenever $z\in A_j$, $j=1\cdots n$, while $\K(\cdot,f_n)=1$ otherwise. As before we also use $|D\phi_A(z)|\leq C\alpha $ when $z\in A$ (and $|D\phi_A(z)|\leq 1$ at all other points) to obtain that $|Df_n|\leq C \alpha_j$ on $A_j$, $j=1\cdots n$, and $|Df_n|\leq 1$ otherwise. Using the exact same arguments as before we see that for the limit
		$
		\overline{f}=\lim_{n\to\infty}f_n
		$
		to be a homeomorphism of integrable distortion it is enough to check that
		\begin{equation}\label{dfinl1general}
			\sum_{n=1}^\infty |A_n|\,\alpha_n<+\infty
		\end{equation}
		and
		\begin{equation}\label{kfinl1general}
			\sum_{n=1}^\infty |A_n|\,\frac{\alpha_n^{2}}{q_n}<+\infty.
		\end{equation}
		Note that as in the case of $f,$ \eqref{kfinl1general} implies \eqref{dfinl1general} when $q_n< \alpha_n$ and so $\alpha_n$ and $q_n$ need to be chosen such that \eqref{kfinl1general} is satisfied as well as the purpose of $\overline{f}$ to be optimal for Corollary \ref{main 1}. 
		Thus we choose 
		\begin{equation}\label{new alpha}
			\alpha_n= \frac{h(r_n)}{r_n}\, \left(\log\frac{1}{r_n}\right)^{1/2}\hspace{2cm}q_n=\log\frac1{r_n},
		\end{equation}
		where $h$ is any gauge function such that $h(r)\to0$ as $r\to0$ and the condition $q_n< \alpha_n$ is satisfied. Indeed, with these choices \eqref{kfinl1general} becomes
		$$\sum_n\left(h(r_n)\right)^{2}<+\infty$$
		which, as before, may always be satisfied by choosing small enough $r_n$. Having \eqref{kfinl1general} fulfilled, our map $\overline{f}$ is a mapping of finite distortion with $\K(\cdot, \overline{f})\in L^1_{loc}$. Furthermore, since we can bound spiraling from below by $\alpha_n$ at the points $r_n$ using the same estimate \eqref{LauriADD} as before, the resulting map $\overline{f}$ attains the optimal rotational behavior stated at Corollary \ref{main 1} modulo the gauge function $h$ which can be chosen to converge to $0$ as slowly as desired.
		Therefore, Theorem \ref{submain 1} will be proven once we show that $\overline{f}$ is H\"older from below. \\
		\\
		To this end, we first observe that the composition of $z\mapsto ze^{i\alpha\log|z|}$ followed by $z\mapsto z|z|^{q-1}$ is precisely $z\mapsto z|z|^{q-1}e^{i\alpha\log|z|}$. This observation suggests us to decompose $\overline{f}=g\circ f$, where $f$ is essentially the first example in this section (with slightly different choices for the constants $\alpha_n$) and $g$ is constructed by building blocks \eqref{buildingblock} with $\alpha=0$ at each step. Morally, $f$ leaves fixed all circles centered at $0$ and only rotates inside the annuli $A_n$, while $g$ conveniently stretches each $A_n$. The H\"older nature of $f^{-1}$ has already been proven when $\alpha_n=\frac{h(r_n)}{r_n}$. We need to show that our map $f^{-1}$ is still H\"older continuous with our new choices for $\alpha_n$, which we can estimate by
		\begin{equation}\label{new estimate}
			\alpha_n= \frac{h(r_n)}{r_n}\, \left(\log\frac{1}{r_n}\right)^{1/2}\leq h(r_n)\,r_n^{-1/(1- \epsilon)}
		\end{equation}
		for an arbitrary $\epsilon>0$ and small enough $r_n$. This can be done by exactly the same proof  as before once we check that $f$ is H\"older from below inside the annuli $A_n.$ To this end, let us consider two points $x,y\in A_n$ and note that $f(x)=f_{n}(x)$ and $f(y)=f_{n}(y)$. Since $r_n>C|x-y|,$ using the estimate \eqref {new estimate} gives
		$$\aligned
		|f(x)-f(y)|=|f_n(x)-f_n(y)| &\gtrsim \frac{r_n^{1/(1- \epsilon)}}{h(r_n)}|x-y|
		\geq \frac{C}{h(r_n)}|x-y|^{1+ \frac{1}{1- \epsilon}}
		\geq C|x-y|^{2+\epsilon}
		\endaligned$$where we are using the bilipschitz nature of $f_n$ in $A_n.$ 
		Therefore, in order to prove Theorem \ref{submain 1} it remains to prove that $g$ is H\"older from below. To this end, given any two points $x,y\in B(0,1)$, we can without loss of generality assume that $|y|\geq|x|$ and let $w$ be the point for which $|w|=|x|$ and $\arg(w)=\arg(y)$. Now, as $g$ is a radial stretching map, it follows that
		$$\aligned
		|g(x)-g(y)|\geq\max\{|g(x)-g(w)|,|g(y)-g(w)|\}.
		\endaligned$$
		Moreover, 
		$
		\max\{|x-w|,|y-w|\}\geq\frac{1}{2}|x-y|.
		$
		Therefore, it is enough to show that both $|g(x)-g(w)|$ and $|g(y)-g(w)|$ satisfy H\"older bounds from below. If $x=0,$ then $w=0$ and we have only $|g(y)-g(w)|$. Let us first check the term $|g(x)-g(w)|$. Since $g$ maps radially circles centered at the origin to similar circles we see that $|g(x)-g(w)|$ gets contracted the same amount as the modulus $|g(x)|$ is contracted under $g$. Now we must consider two possibilities, either $x,w \in A_n$ or $x,w \in D_n$ for some $n$. Let us first assume $x,w\in A_n=B(0,R_n)\setminus B(0,r_n)$ for some $n$. Then following exactly the same argument as in the proof of Theorem 3 in \cite{CHS} (page no. 19), one can obtain the estimate
		%Here we recall the ansatz \eqref{Additional assumption Lauri} on $r_n$. Then
		
		%$$\aligned
		%|g(x)|=\left(\frac{1}{e}\right)^{q_{n-1}+...+q_1-(n-1)}\cdot |x|\,\left(\frac{|x|}{R_n}\right)^{q_n-1}
		%\geq r_n \cdot |x|\,\left(\frac{|x|}{R_n}\right)^{q_n-1}
		%\geq r_n \cdot |x|\,\left(\frac1e\right)^{q_n-1} =e\cdot r_{n}^{1+ \beta}\cdot|x|
		%\endaligned$$ for any $x\in A_n$, where in the last step we use \eqref{new alpha}. Therefore,
		$$\aligned
		|g(x)-g(w)|\geq C\cdot|x-w|^3,
		\endaligned$$
		for some fixed constant $C>0$ whenever $x,w \in A_n$ or $x,w \in D_n$. %Next, let $x,w\in D_n=B(0,r_n)\setminus B(0,R_{n+1})$ for some $n$.  Using \eqref{Additional assumption Lauri} we get
		%$$\aligned
		%|g(x)|&\geq c\left(\frac{1}{e}\right)^{q_{n-1}+...+q_{1}-(n-1)}\cdot r_{n}^{\beta}\cdot|x|
		% \geq c\cdot r_{n}^{1+ \beta}\cdot|x|.
		%\endaligned$$
		%Thus we can use similar argument as in the previous case to estimate
		%$$\aligned
		%|g(x)-g(w)|&\geq c\cdot r_{n}^{1 +\beta}\cdot|x-w|
		%\geq c\cdot|x-w|^{2+ \beta}\\
		%\endaligned$$since $|x-w|<c\cdot r_n$ for some fixed constant $c>0$ when $x,w\in D_n$.\\  
		\\
		Since the set $\D\setminus\{0\}$ is partitioned by separated annuli $A_n$ and $D_n$ it is clear that $|g(x)-g(w)|$ satisfies H\"older estimates from below. \\
		Finally, let us prove the H\"older estimates from below for the term $|g(y)-g(w)|$. As the mapping $g$ is radial, we can  assume that $y$ and $w$ are real.
		Now, we estimate the differential of $g$ from below using similar technique as in the proof of Theorem 3 in \cite{CHS} (page no. 20),
		%We intend to use the Fundamental Theorem of Calculus, and thus have to estimate the differential from below.  Using \eqref{Additional assumption Lauri}, as well as the facts that $q_n>1$ and $R_n=er_n$, we can estimate for any real number $t\in [r_n,R_n]$ that
		%$$\aligned
		%g'(t)=\left(\frac{1}{e}\right)^{q_{n-1}+...+q_1-(n-1)}\cdot   q_n\cdot \left(\frac{t}{R_n}\right)^{q_{n}-1}
		%\geq r_n\, q_n\cdot \left(\frac{r_n}{R_n}\right)^{q_{n}-1}   
		%= e\,q_n\,r_{n}^{1+ \beta} 
		%\geq c \beta\cdot t^{1+ \beta}\,\log\left(e+ %\frac1t \right).
		%\endaligned$$ 
		%Next, if $t\in [R_{n+1},r_n]$, we have
		$$\aligned
		g'(t) \geq c \cdot t^2
		\endaligned$$
		for some constant $c>0,$ if $t\in [r_n,R_n]$ or $t\in [R_{n+1},r_n]$.
		%$$\aligned
		%g'(t)=\left(\frac{1}{e}\right)^{q_{n-1}+...+q_{1}-(n-1)}\cdot\left(\frac{1}{e}\right)^{q_{n}-1}\geq e\cdot r_{n}^{1+ \beta}\geq c\cdot t^{1+ \beta}
		%\endaligned$$
		Thus, as above, since $(0,1)$ is partitioned by the intervals $[r_n,R_n]$, $[R_{n+1},r_n]$ and $[R_1,1)$, one gets that   
		$$\aligned
		g'(t)\geq c\cdot t^2
		\endaligned$$ for every $t\in(0,1)$. Now, we use the fundamental theorem of calculus to get
		$$\aligned
		|g(y)-g(w)|=\int_{w}^{y}g'(t)dt
		\geq\int_{w}^{y}c\cdot t^2 dt
		=C\left(y^{3}-w^{3}\right)
		&\geq C|y-w|^{3}.
		\endaligned$$ This proves that the second term is H\"older from below as well, which in turn proves that $g$ is H\"older from below with exponent $3$. This finishes the proof of Theorem \ref{submain 1}.
	\end{proof}
	
	%\noindent\textbf{Conflict of Interest} : The authors declare that they have no conflicts of interest. \\ \\
	%\textbf{Data availability statement} : The results presented in this study are openly available in Arxiv at arXiv:2110.12809. However, the results in this manuscript are a bit modified from those in the arxiv version.
	
	\vspace{1cm}
	Lauri Hitruhin\\
	Department of Mathematics and Statistics, University of Helsinki, PO Box 68, 00014, Helsinki, Finland \\
	Email address: lauri.hitruhin@gmail.com \\ \\
	
	\noindent Banhirup Sengupta\\
	Centre for Applicable Mathematics, Tata Institute of Fundamental Research, PO Box 6503, GKVK Post
	Office, Bangalore 560065, India \\
	Email address: sengupta25@tifrbng.res.in

\end{document}